\documentclass[12pt]{article}
\usepackage{graphicx}
\usepackage{color}
\usepackage{color}
\catcode`\@=11 \@addtoreset{equation}{section}

\catcode`\@=12

\newtheorem{Theorem}{Theorem}[section]
\newtheorem{Proposition}{Proposition}[section]
\newtheorem{Lemma}{Lemma}[section]
\newtheorem{Corollary}{Corollary}[section]

\newcommand{\bTheorem}[1]{
\begin{Theorem} \label{T#1} }
\newcommand{\eT}{\end{Theorem}}

\newcommand{\bProposition}[1]{
\begin{Proposition} \label{P#1}}
\newcommand{\eP}{\end{Proposition}}

\newcommand{\bLemma}[1]{
\begin{Lemma} \label{L#1} }
\newcommand{\eL}{\end{Lemma}}

\newcommand{\bCorollary}[1]{
\begin{Corollary} \label{C#1} }
\newcommand{\eC}{\end{Corollary}}

\newcommand{\bFormula}[1]{
\begin{equation} \label{#1}}
\newcommand{\eF}{\end{equation}}

\newcommand{\Ov}[1]{\overline{#1}}

\newcommand{\DC}{C^\infty_c}

\newcommand{\vr}{\varrho}
\newcommand{\vre}{\vr_\ep}

\newcommand{\vue}{\vu_\ep}

\newcommand{\vu}{\vc{u}}
\newcommand{\vc}[1]{{\bf #1}}

\newcommand{\Div}{{\rm div}_x}
\newcommand{\Grad}{\nabla_x}

\newcommand{\tn}[1]{\mbox {\F #1}}
\newcommand{\dx}{{\rm d} {x}}
\newcommand{\dt}{{\rm d} t }

\newcommand{\intO}[1]{\int_{\Omega} #1 \ \dx}

\newcommand{\ep}{\varepsilon}

\font\F=msbm10 scaled 1000

\definecolor{Cgrey}{rgb}{0.85,0.85,0.85}
\definecolor{Cblue}{rgb}{0.50,0.85,0.85}
\definecolor{Cred}{rgb}{1,0,0}
\definecolor{fancy}{rgb}{0.10,0.85,0.10}

\newcommand\Cbox[2]{%
    \newbox\contentbox%
    \newbox\bkgdbox%
    \setbox\contentbox\hbox to \hsize{%
        \vtop{
            \kern\columnsep
            \hbox to \hsize{%
                \kern\columnsep%
                \advance\hsize by -2\columnsep%
                \setlength{\textwidth}{\hsize}%
                \vbox{
                    \parskip=\baselineskip
                    \parindent=0bp
                    #2
                }%
                \kern\columnsep%
            }%
            \kern\columnsep%
        }%
    }%
    \setbox\bkgdbox\vbox{
        \color{#1}
        \hrule width  \wd\contentbox %
               height \ht\contentbox %
               depth  \dp\contentbox
        \color{black}
    }%
    \wd\bkgdbox=0bp%
    \vbox{\hbox to \hsize{\box\bkgdbox\box\contentbox}}%
    \vskip\baselineskip%
}

\date{}

\begin{document}


\title{Scale interactions in compressible rotating fluids}

\author{Eduard Feireisl\thanks{Eduard Feireisl acknowledges the support of the project LL1202 in the
programme ERC-CZ funded
by the Ministry of Education, Youth and Sports of the Czech Republic.} \and
Anton\' \i n Novotn\' y }

\maketitle

\bigskip

\centerline{Institute of Mathematics of the Academy of Sciences of the Czech Republic}

\centerline{\v Zitn\' a 25, 115 67 Praha 1, Czech Republic}

\centerline{Charles University in Prague, Faculty of Mathematics and Physics, Mathematical Institute}

\centerline{Sokolovsk\' a 83, 186 75 Praha 8,
Czech Republic}

\centerline{and}

\centerline{IMATH, EA 2134, Universit\' e du Sud Toulon-Var
BP 20132, 83957 La Garde, France}

\begin{abstract}

We study a triple singular limit for the scaled barotropic Navier-Stokes system modeling the motion of a rotating,
compressible, and viscous fluid, where the Mach and Rossby numbers are proportional to a small parameter $\ep$, while
the Reynolds number becomes infinite for $\ep \to 0$. If the fluid is confined to an infinite slab
bounded above and below by two parallel planes, the limit behavior is identified as a purely horizontal motion of an
incompressible inviscid fluid, the evolution of which is described by an analogue of the  Euler system.

\end{abstract}

\tableofcontents

\section{Introduction}
\label{i}

There have been several types of singular limits of the Navier-Stokes system studied in the recent literature, among
them a few devoted to the effect of rotation, see Chemin et al. \cite{CDGG}. In this paper, we consider the scaled
\emph{compressible Navier-Stokes system} describing the time evolution
of the \emph{density} $\vr = \vr(t,x)$ and the \emph{velocity} $\vu = \vu(t,x)$ of a compressible viscous
and rotating fluid:
\bFormula{i1}
\partial_t \vr + \Div( \vr \vu ) = 0,
\eF
\bFormula{i2}
\partial_t (\vr \vu) + \Div (\vr \vu \otimes \vu) + \frac{1}{{\rm Ro}} \vr ( \vc{ \omega} \times \vu)
+ \frac{1}{{\rm Ma}^2} \Grad p(\vr) = \frac{1}{{\rm Re}} \Div \tn{S} + \vr \Grad G,
\eF
where $\tn{S}$ is the viscous stress, here given by \emph{Newton's rheological law},
\bFormula{i3}
\tn{S} =
\left[ \mu \left( \Grad \vu + \Grad^t \vu - \frac{2}{3} \Div \vu \tn{I} \right) + \eta \Div \vu \right],
\eF
$p = p(\vr)$ is the pressure, $\mu > 0$ and $\eta \geq 0$ are the
viscosity coefficients, $\vc{\omega} = [0,0,1]$ is the axis of rotation, and $\Grad G$ represents a conservative force imposed on the system, say, by the
gravitational potential $G$ of an object placed \emph{outside} the fluid domain, see the survey of Klein \cite{Klein}.

The scaled system contains several
characteristic numbers:

\begin{itemize}
\item
{\rm Ro} \dotfill \mbox{Rossby number}

\item
{\rm Ma} \dotfill \mbox{Mach number}

\item
{\rm Re} \dotfill \mbox{Reynolds number}

\end{itemize}

The following are examples of singular limits considered in numerous studies:

\bigskip

\begin{itemize}

\item
{\textsc{The low Mach number limit:}} The Mach number is the ratio of the characteristic speed of the fluid divided on the speed of sound. In the low Mach number limit, the fluid flow becomes \emph{incompressible}, the density distribution is constant and the velocity field solenoidal, see Ebin \cite{EB1}, Klainerman and Majda \cite{KM1},
Lions and Masmoudi \cite{LIMA1}, Masmoudi \cite{MAS1}, among others.

\item
{\textsc{The low Rossby number limit:}} Low Rossby number corresponds to fast rotation.
As observed by many authors, the highly rotating fluids become \emph{planar} (two-dimensional). Accordingly, the fast rotation has a regularizing effect, see Babin et al.
\cite{BaMaNi2}, \cite{BaMaNi1},  Chemin et al. \cite{CDGG}.

\item
{\textsc{The high Reynolds number limit:}} In the high Reynolds number limit, the viscosity of the fluid becomes negligible. Consequently, solutions of the Navier-Stokes system tend to the solutions of the Euler system, see
Clopeau, Mikeli{\' c}, Robert \cite{CloMikRob}, Masmoudi \cite{MAS6}, \cite{MAS7}, \cite{MAS1}, Swann
\cite{Swa}, among others. The inviscid limits include the difficulties related to the boundary behavior of the fluid and a proper choice of boundary conditions, see Kato \cite{Kato}, Kelliher \cite{Kelli1}, Sammartino and Caflisch
\cite{SamCaf1}, \cite{SamCaf2}, Temam and Wang \cite{TemWan1}, \cite{TemWan2}.

\end{itemize}

\bigskip

The effects described above may act \emph{simultaneously}. The incompressible inviscid limit was investigated by Masmoudi \cite{MAS7}, for viscous rotating fluids see Masmoudi \cite{MAS8}, Ngo \cite{Ngo}, the compressible rotating fluids were discussed in \cite{FeGaGVNo}, \cite{FeGaNo}. In this paper, we address the problem of the triple limit for
${\rm Ma} = {\rm Ro} = \ep \to 0$, while ${\rm Re} = {\rm Re}(\ep) \to \infty$ as $\ep \to 0$. In agreement with the previous discussion, the fluid flow is expected to become {\bf (i)} incompressible, {\bf (ii)} planar (2D), and {\bf (iii)} inviscid and as such described by a variant of the 2D incompressible Euler system that is known to possess global-in-time solutions for any regular initial data. Note that the action of volume forces in the momentum equation (\ref{i2})
is represented solely by the potential $G$, notably
the effect of the \emph{centrifugal} force is neglected. This is a standard simplification adopted, for instance, in
models of atmosphere or astrophysics, see Jones et al. \cite{JoRo}, \cite{JoRoGa} , Klein \cite{KL2}. On the other hand, although the centrifugal force is
counterbalanced by the gravity in many real-world applications (see Durran \cite{Dur1}), it is proportional to
$1/\ep^2$ under the present scaling, and its far-field impact may change the limit problem dramatically, see \cite{FeGaGVNo}.

Similarly to \cite{FeGaGVNo}, \cite{FeGaNo}, we consider the problem (\ref{i1} - \ref{i3}) in an infinite slab
$\Omega = R^2 \times (0,1)$, denoting the \emph{horizontal components} of a vector field $\vc{v}$ by
$\vc{v}_h$, $\vc{v} = [\vc{v}_h, v_3]$, where the velocity field $\vu$ satisfies the complete slip boundary
conditions
\bFormula{i4}
\vc{u} \cdot \vc{n} = u_3|_{\partial \Omega} = 0, \ [\tn{S} (\Grad \vu) \vc{n}]_{\rm tan} =
[S_{2,3}, - S_{1,3}, 0]|_{\partial \Omega} = 0,
\eF
where $\vc{n} = [0,0, \pm 1]$ is the outer normal vector. Such a choice of boundary behavior prevents the flow from creating a viscous boundary layer - the up to now unsurmountable difficulty of the inviscid limits, see Kato \cite{Kato}, Temam and Wang \cite{TemWan2}. As a matter of fact, replacing
(\ref{i4}) by the more standard no-slip boundary condition would drive the fluid to the trivial state $\vu = 0$ in the asymptotic limit unless we impose anisotropic viscosity, see Bresch, Desjardins and Gerard-Varet \cite{BrDeGV},
Chemin et al. \cite{CDGG}.

Our approach is based on the \emph{relative entropy inequality} (cf. \cite{FeJiNo}, Masmoudi \cite{MAS7}, Wang and Jiang \cite{WanJia}) applied in the framework of weak solutions to the Navier-Stokes system (\ref{i1} - \ref{i3}).
We consider the \emph{ill-prepared} initial data:
\bFormula{i5}
\vr(0, \cdot) = \vr_{0, \ep} = \Ov{\vr} + \ep \vr^{(1)}_{0, \ep}, \ \vr^{(1)}_{0,\ep} \to \vr^{(1)}_0 \ \mbox{in} \ L^2(\Omega),
\eF
\bFormula{i6}
\vu(0, \cdot) = \vu_{0, \ep} \to \vu_0 \ \mbox{in} \ L^2(\Omega; R^3),
\eF
where $\Ov{\vr} > 0$ is the anticipated constant limit density enforced by the incompressible limit. Accordingly,
the resulting problem is supplemented by the far field conditions
\bFormula{i7}
\vr \to \Ov{\vr}, \ \vu \to 0 \ \mbox{as} \ |x_h| \to \infty.
\eF

Supposing we already know that, in some sense,
\[
\vr^{(1)}_\ep = \frac{\vre - \Ov{\vr}}{\ep} \to q, \ \vue \to \vc{v}
\]
we may (formally) check that $q = q(x_h)$, $\vc{v} = [\vc{v}_h (x_h),0]$ are interrelated through the \emph{diagnostic equation}
\bFormula{i8}
\vc{\omega} \times \vc{v} + \frac{p'(\Ov{\vr} )}{\Ov{\vr}} \Grad q = 0,
\eF
and $q$ satisfies
\bFormula{i9}
\partial_t \left( \Delta_h q - \frac{1}{p'(\Ov{\vr})} q \right) +
\nabla_h^\perp q \cdot \nabla_h \left( \Delta_h q -
\frac{1}{p'(\Ov{\vr})} q \right) = 0. \eF Here and hereafter, the
subscript $h$ indicates the restriction of the standard
differential operators to the horizontal variables, for instance,
$\nabla_h f = [\partial_{x_1} f , \partial_{x_2} f]$,\ ${\rm
div}_h \vc{v} = \partial_{x_1} h_1 + \partial_{x_2} h_2$,
$\Delta_h = {\rm div}_h \nabla_h$, etc.

Note that
\[
\nabla_h^\perp q = \frac{\Ov{\vr}}{p'(\Ov{\vr})} \vc{v}, \ \Delta_h q = {\rm curl}_h \vc{v}_h;
\]
whence $r$ can be viewed as a kind of \emph{stream function}, while the system (\ref{i8}), (\ref{i9}) possesses
the same structure as the 2D Euler equations. In particular, we expect the solutions of (\ref{i8}), (\ref{i9}) to
be as regular as the initial data and to exist globally in time. Equation (\ref{i9}) arises in the
theory of quasi-geostrophic flows, see Zeitlin \cite[Chapters 1,2]{Zeit}.

One of the major stumbling blocks in the analysis of the singular limit is the presence of rapidly oscillating
Rossby-acoustic waves. Their behavior is described by means of a hyperbolic system
\bFormula{i10}
\ep \partial_t s + \Ov{\vr} \Div \vc{V} = 0,
\eF
\bFormula{i11}
\ep \partial_t \vc{V} + \left( \omega \times \vc{V} + \frac{p'(\Ov{\vr})}{\Ov{\vr} } \Grad s \right) = 0,
\eF
where
\[
s = s_\ep \approx \frac{\vre - \Ov{\vr}}{\ep}, \ \vc{V} \approx \vc{V}_\ep = \vue.
\]

The bulk of the paper is devoted to the dispersive estimates for the problem (\ref{i10}), (\ref{i11}). In particular, we use the recent results of Guo, Peng, and Wang \cite{GuPeWa} on the asymptotic behavior of the abstract group
of operators
\[
t \mapsto \exp \left[ {\rm i} t \Phi \left( \sqrt{ - \Delta } \right) \right],
\]
where $\Phi$ is a function with specific properties. In particular, we establish
$L^1 - L^\infty$ decay estimates for the solutions of (\ref{i10}), (\ref{i11}) in the frequency domain bounded away from zero.

The paper is organized as follows. In Section \ref{m}, we introduce the standard definition of \emph{finite-energy weak solutions} to the scaled system (\ref{i1} - \ref{i6}) and formulate the main result. Section \ref{r} contains the relative entropy inequality, together with the uniform bounds on the family of solutions of the scaled system.
Section \ref{d} represents the bulk of the paper. Using the abstract result of Guo et al. \cite{GuPeWa} we establish
the $L^1-L^\infty$ estimates for the acoustic-Rossby waves.
Such a result for the system (\ref{i10}), (\ref{i11}) may be of independent interest and represents an analogue of the standard Strichartz estimates for the wave and Schroedinger equations. In particular, we extend the smoothing estimates established in \cite{FeGaNo} and
obtain the necessary tool to attack the inviscid limit in Section \ref{c}.

\section{Preliminaries, main results}
\label{m}

In order to fix ideas and to simplify presentation, we suppose, without loss of generality, that $\Ov{\vr} = 1$.
In addition, we assume that the pressure $p \in C [0,\infty) \cap C^3 (0, \infty)$ satisfies
\bFormula{m1}
p(0) = 0, \ p'(\vr) > 0 \ \mbox{for all}\ \vr > 0,\ \lim_{\vr \to \infty} \frac{p'(\vr)}{\vr^{\gamma - 1}} =
p_\infty > 0 \ \mbox{for a certain}\ \gamma > \frac{3}{2}.
\eF
Morever, again for the sake of simplicity, we suppose
\bFormula{m2}
p'(\Ov{\vr}) = 1.
\eF

Finally, given our choice of the complete slip boundary conditions (\ref{i4}), it is convenient to
replace the set $\Omega = R^2 \times [0,1]$ by
\[
\Omega = R^2 \times [-1,1]|_{\{ -1, 1 \}},
\]
meaning we suppose that all quantities are 2-\emph{periodic} with respect to the vertical variable $x_3$. Moreover, in accordance with (\ref{i4}), we assume that
\bFormula{m2a}
\vr(t,x_h, - x_3) = \vr(t, x_h, x_3),\
\vc{u}_h (t,x, - x_3) =  \vc{u}_h (t,x, x_3), \
u_3 (t,x, - x_3) = - u_3 (t,x, - x_3),
\eF
and
\bFormula{m2bb}
G(x_h, -x_3) = G(x_h, x_3)
\eF
for all $t \in (0,T)$, $x_h \in R^2$, $x_3 \in [-1,1]|_{\{
-1,1 \} }$.
Such a formulation, completely equivalent to (\ref{i4}), was proposed by Ebin \cite{EB}.

Setting ${\rm Ma} = {\rm Ro} = \ep$, $\mu = \mu_\ep \searrow 0$ we say that
$\vr$, $\vu$ is a \emph{finite energy weak solution} to the scaled Navier-Stokes system (\ref{i1} - \ref{i7}) if:

\begin{itemize}
\item The density $\vr$ is a non-negative function such that
\[
(\vr - 1) \in L^\infty(0,T; (L^2 + L^{\gamma})(\Omega));
\]
the velocity $\vu$ belongs to the space $L^2(0,T; W^{1,2}(\Omega))$. Moreover, in accordance with our convention,
both $\vr$ and $\vu$ satisfy the symmetry condition (\ref{m2a}).

\item The equation of continuity (\ref{i1}) holds in the weak sense:
\bFormula{m3}
\int_0^T \intO{ \left( \vr \partial_t \varphi + \vr \vu \cdot \Grad \varphi \right) } \ \dt =
- \intO{ \vr_{0,\ep} \varphi(0, \cdot) }
\eF
for any $\varphi \in \DC([0,T) \times {\Omega})$.

\item Similarly, the momentum equation is replaced by a family of integral identities
\bFormula{m4}
\int_0^T \intO{ \left( \vr \vu \cdot \partial_t \varphi + \vr \vu \otimes \vu : \Grad \varphi -
\frac{1}{\ep} \vr (\vc{\omega} \times \vu) \cdot \varphi + \frac{1}{\ep^2} p(\vr) \Div \varphi \right) }
\ \dt
\eF
\[
=\int_0^T \intO{ \tn{S}_\ep (\Grad \vu) : \Grad \varphi } \ \dt - \int_0^T \intO{ \vr \Grad G \cdot \varphi} \ \dt -  \intO{ \vr_{0,\ep} \vu_{0,\ep} \cdot \varphi(0, \cdot)
}
\]
for any $\varphi \in \DC([0,T) \times {\Omega})$,
with
\bFormula{m5}
\tn{S}_\ep (\Grad \vu) = \mu_\ep \left( \Grad \vu + \Grad^t \vu - \frac{2}{3} \Div \vu \tn{I} \right),\
\mu_\ep \searrow 0.
\eF

\item
The energy inequality
\bFormula{m6}
\intO{ \left[ \frac{1}{2} \vr |\vu|^2 + \frac{1}{\ep^2} \left( H(\vr) - H'(1) (\vr - 1) - H(1) \right) \right] (\tau, \cdot) }
+ \int_0^\tau \intO{ \tn{S}_\ep (\Grad \vu) : \Grad \vu } \ \dt
\eF
\[
\leq \intO{ \left[ \frac{1}{2} \vr_{0,\ep} |\vu_{0,\ep} |^2 + \frac{1}{\ep^2} \left( H(\vr_{0,\ep}) - H'(1) (\vr_{0,\ep} - 1) - H(1) \right) \right] }
+ \int_0^\tau \intO{ \vr \Grad G \cdot \vu } \ \dt
\]
holds for a.a. $\tau \in [0,T]$, where we have set
\bFormula{m6a}
H(\vr) = \vr \int_1^\vr \frac{p(z)}{z^2} \ {\rm d}z.
\eF
\end{itemize}

\subsection{Limit system}

Under the convention (\ref{m1}), (\ref{m2}), the expected limit problem reads
\bFormula{m7}
\vc{\omega} \times \vc{v} + \Grad q = 0, \ \vc{v} = [\vc{v}_h(x_h),0] ,\ q = q(x_h),
\eF
\bFormula{m8}
\partial_t \Big( \Delta_h q - q \Big) + \vc{v}_h \cdot \nabla_h (\Delta_h q ) = 0,
\eF
supplemented with the initial condition
\bFormula{m9}
q(0,\cdot) = q_0.
\eF

Note that (\ref{m8}) can be written as
\bFormula{m10}
\partial_t \Big( \Delta_h q - q \Big) + \vc{v}_h \cdot \nabla_h \Big( \Delta_h q - q \Big) = 0;
\eF
whence the problem enjoys strong similarity with the standard Euler system. In particular, we may use the abstract theory
of Oliver \cite[Theorem 3] {Oli} to obtain the following result:

\bProposition{m1}
Suppose that
\[
q_0 \in W^{m,2}(R^2) \ \mbox{for}\ m \geq 4.
\]

Then the problem (\ref{m9}), (\ref{m10}) admits a solution $q$, unique in the class
\[
q \in C([0,T]; W^{m,2}(R^2) \cap C^1([0,T]; W^{m-1,2}(R^2)).
\]
\eP

It is worth noting that, similarly to the $2D$-Euler system, the solution $q$ can be constructed \emph{globally} in time.

\subsection{Main result}

Having collected all the necessary preliminary material we are in a position to state the main result of the present paper.

\bTheorem{m1}
Let the pressure $p$ satisfy the hypotheses (\ref{m1}), (\ref{m2}). Suppose that the initial data $\vr_{0,\ep}$,
$\vu_{0,\ep}$ belong to the symmetry class (\ref{m2a}) and are given through (\ref{i5}), (\ref{i6}), where
\bFormula{m11-}
\{ \vr^{(1)}_{0,\ep} \}_{\ep > 0} \ \mbox{bounded in} \ L^2 \cap L^\infty (\Omega),\
\vr^{(1)}_{0,\ep} \to \vr^{(1)}_0 \ \mbox{in}\ L^2(\Omega),
\eF
\bFormula{m11--}
\{ \vu_{0,\ep} \}_{\ep > 0} \ \mbox{bounded in} \ L^2 (\Omega;R^3),\
\vu_{0,\ep} \to \vu_0 \ \mbox{in}\ L^2(\Omega;R^3),
\eF
where
\bFormula{m11a}
\vr^{(1)}_0 \ \in W^{m-1,2}(\Omega),\
\vu_0 \in W^{m,2}(\Omega;R^3), \ m \geq 3 .
\eF
Let
\bFormula{hypp}
\Grad G \in L^\infty \cap L^r (\Omega;R^3) \ \mbox{for a certain} \ 1 \leq r < 2,
\eF
and satisfy (\ref{m2bb}).

Furthermore, let
$q_0 = q_0(x_h)$ be the unique solution of the elliptic problem
\bFormula{m11}
- \Delta_h q_0 + q_0 = \int_0^1 {\rm curl}_h [\vu_0]_h \ {\rm d}x_3 + \int_0^1 \vr^{(1)}_0 \ {\rm d}x_3
\ \mbox{in}\ W^{1,2}(R^2).
\eF

Finally,
let $[\vre, \vue]$ be a weak solution of the scaled Navier-Stokes system (\ref{i1} - \ref{i7}) in the sense specified above.

Then
\[
\frac{\vre - 1}{\ep} \to q  \left\{ \begin{array}{l} \mbox{weakly-(*) in}\ L^\infty(0,T; (L^2 + L^r)(\Omega))
\\ \\ \mbox{(strongly) in} \ L^1_{\rm loc}((0,T) \times {\Omega}),
\end{array} \right.
\]
\[
\sqrt{\vre} \vue \to \vc{v} \left\{ \begin{array}{l} \mbox{weakly-(*) in} \ L^\infty(0,T;
L^2  (\Omega; R^3)), \\ \\ \mbox{strongly in}\ L^1_{\rm loc}((0,T) \times {\Omega};R^3),
\end{array} \right.
\]
where $[q, \vc{v}]$ is the (unique) solution of the problem (\ref{m7} - \ref{m9}).

\eT

The rest of the paper is devoted to the proof of Theorem \ref{Tm1}.

\section{Relative entropy, uniform bounds}
\label{r}

We start by introducing the relative entropy functional for the compressible Navier-Stokes system identified in
\cite{FeJiNo}, \cite{FeNoSun}, Germain \cite{Ger}. Set
\bFormula{rr1}
\mathcal{E}_\ep \left( \vr, \vu \Big| r, \vc{U} \right) =
\intO{  \left[ \frac{1}{2} \vre |\vue - \vc{U}|^2 + \frac{1}{\ep^2} \Big( H(\vre) - H'(r)(\vre - r) - H(r) \Big)(\tau, \cdot) \right]},
\eF
where the function was defined in (\ref{m6a}).

\subsection{Relative entropy inequality}

It can be shown that \emph{any} finite energy weak solution $[\vre, \vue]$ of the Navier-Stokes system (\ref{i1} - \ref{i7})
satisfies the \emph{relative entropy inequality}:
\bFormula{r1}
\mathcal{E}_\ep \left( \vre, \vue \ \Big| \ r, \vc{U} \right) (\tau)
+ \int_0^\tau \intO{ \Big( \tn{S}_\ep (\Grad \vue) - \tn{S}_\ep (\Grad \vc{U}) \Big) : \Big( \Grad \vue - \Grad \vc{U} \Big) } \ \dt \leq
\eF
\[
\mathcal{E}_\ep \left( \vr_{0,\ep}, \vu_{0,\ep} \ \Big| \ r(0,\cdot) , \vc{U}(0,\cdot) \right)
\]
\[
+ \int_0^\tau \intO{  \vre \left( \partial_t \vc{U} + \vue \cdot \Grad \vc{U} \right) \cdot \left( \vc{U} - \vue \right) } \ \dt
\]
\[
+ \int_0^\tau \intO{ \tn{S}_\ep (\Grad \vc{U}) : \Grad (\vc{U} - \vue ) } \ \dt + \frac{1}{\ep} \int_0^\tau \intO{
\vre (\vc{\omega} \times \vue ) \cdot (\vc{U} - \vue) } \ \dt
\]
\[
+ \frac{1}{\ep^2} \int_0^\tau \intO{ \Big[ (r - \vre) \partial_t H'(r) + \Grad H'(r) \cdot (r \vc{U} - \vre \vue ) \Big] } \ \dt
- \frac{1}{\ep^2} \int_0^\tau \intO{ \Div \vc{U} \Big( p(\vre) - p(r) \Big) } \ \dt
\]
\[
- \int_0^\tau \intO{ \vre \Grad G \cdot (\vc{U} - \vue) } \ \dt
\]
for all (smooth) functions $r$, $\vc{U}$ such that
\bFormula{r1a}
r > 0, \ (r-1) \in \DC([0,T] \times \Ov{\Omega}),\
\vc{U} \in \DC([0,T] \times {\Omega}),
\eF
see \cite{FeJiNo}. Clearly, the class of admissible ``test'' functions (\ref{r1a}) can be considerably extended by means of a density argument. Note that the relative entropy inequality (\ref{r1}) reduces to the energy inequality
(\ref{m6}) provided we take $r = 1$, $\vc{U} = 0$.

\subsection{Uniform bounds}

Before deriving the available uniform bounds on the family of solutions $[\vre, \vue]$ it seems convenient to introduce
the essential and residual component of any function $h$:
\[
h = h_{\rm ess} + h_{\rm res},
\]
\[
h_{\rm ess} = \chi(\vre) h,\ \chi \in \DC(0, \infty), \ 0 \leq \chi \leq 1,\
\chi = 1 \ \mbox{on an open interval contaning} \ \Ov{\vr} = 1,
\]
\[
h_{\rm res} = (1 - \chi(\vre)) h,
\]
cf. \cite{FeGaGVNo}, \cite{FeGaNo}.

The uniform bounds are derived from the energy inequality (\ref{m6}) (the relative entropy inquality
(\ref{r1}) with $r = 1$, $\vc{U} = 0$). Since the initial data satisfy the hypotheses
(\ref{m11-}), (\ref{m11--}), the integral on the left-hand side of (\ref{m6}) remains bounded uniformly for
$\ep \to 0$. Accordingly, we get
\bFormula{rr2}
{\rm ess} \sup_{t \in (0,T)} \left\| \sqrt{\vre} \vue \right\|_{L^2(\Omega;R^3)} \leq c,
\eF
\bFormula{rr3}
{\rm ess} \sup_{t \in (0,T)} \left\| \left[ \frac{\vre - 1}{\ep} \right]_{\rm ess} \right\|_{L^2(\Omega)} \leq c,
\eF
\bFormula{rr4}
{\rm ess} \sup_{t \in (0,T)} \left\| [\vre ]_{\rm res} \right\|^\gamma_{L^\gamma(\Omega)}
+ {\rm ess} \sup_{t \in (0,T)} \left\| [ 1 ]_{\rm res}  \right\|^\gamma_{L^\gamma(\Omega)}
\leq \ep^2 c,
\eF
and
\bFormula{rr5}
\mu_{\ep} \int_0^T \intO{ \left| \Grad \vue + \Grad \vue^t - \frac{2}{3} \Div \vue \tn{I} \right|^2 } \ \dt \leq c,
\eF
see \cite[Section 2]{FeGaNo}. We remark that
\[
\left| \intO{ \vre \Grad G \cdot \vue } \right| \leq c \intO{ \left( \vre |\vue|^2 + \left( [\vre]^\gamma_{\rm res} + 1 \right) |\Grad G|^2 \right) },
\]
where, thanks to the hypothesis (\ref{hypp}), the right integral can be ``absorbed'' by means of a Gronwall-type argument.

Finally, by virtue of Korn's inequality, the relation (\ref{rr5}) implies
\bFormula{rr6}
\mu_{\ep} \int_0^T \intO{ |\Grad \vue |^2 } \ \dt \leq c.
\eF

Note that all the above bounds depend, in general, on $T$.

\subsection{Convergence, part I}

It follows immediately from (\ref{rr2} - \ref{rr4}) that
\bFormula{rr7}
\vr^{(1)}_\ep \equiv \frac{ \vre - 1 }{\ep} \to \vr^{(1)}
\ \mbox{weakly-(*) in} \ L^\infty(0,T; (L^2 + L^r) (\Omega)), \ r = \min\{ 2, \gamma \},
\eF
\bFormula{rr8}
\sqrt{\vre} \vue \to \vc{u} \ \mbox{weakly-(*) in}\ L^\infty(0,T; L^2(\Omega;R^3))
\eF
at least for suitable subsequences.
In particular,
\bFormula{rr9}
\vre \to 1 \ \mbox{in} \ L^\infty(0,T; (L^2 + L^r) (\Omega)), \ r = \min\{ 2, \gamma \},
\eF
and
\bFormula{rr10}
\vre \vue \to \vc{u} \ \mbox{weakly-(*) in}\ L^\infty(0,T; (L^2 + L^{2 \gamma/(\gamma + 1)}(\Omega;R^3)).
\eF

Finally, letting $\ep \to 0$ in the equation of continuity (\ref{m3}) we deduce
\bFormula{rr11}
\Div \vc{u} = 0 ,
\eF
and, multiplying the momentum equation by $\ep$, we get the diagnostic equation
\bFormula{rr12}
\vc{\omega} \times \vc{u} + \Grad \vr^{(1)} = 0,
\eF
where both relations are to be understood in the sense of distributions. It is easy to check that (\ref{rr12})
imposes the following restrictions:
\bFormula{rr13}
\vr^{(1)} \ \mbox{independent of}\ x_3, \ \mbox{meaning} \ \vr^{(1)} = \vr^{(1)}(x_h),
\eF
\bFormula{rr14}
\vc{u} = \vc{u}(x_h), \ \Div \vc{u} = {\rm div}_h \vc{u}_h = 0.
\eF
Finally, since $\vc{u}$ belongs to the symmetry class (\ref{m2a}),
\bFormula{rr15}
u_3 = 0, \ \vc{u} = [ \vc{u}_h, 0].
\eF

\section{Dispersive estimates}
\label{d}

As already pointed out in the introduction, the heart of the paper are dispersive estimates for the
acoustic-Rossby waves, the propagation of which is governed by the system
\bFormula{d1}
\ep \partial_t s + \Div \vc{V} = 0,\ s(0, \cdot) = s_0,
\eF
\bFormula{d2}
\ep \partial_t \vc{V} + \vc{\omega} \times \vc{V} + \Grad s = 0, \ \vc{V}(0,\cdot) = \vc{V}_0.
\eF

\subsection{The wave propagator}

Consider the operator
\[
\mathcal{B}: \left[ \begin{array}{c} s \\ \vc{V} \end{array} \right] \mapsto \left[ \begin{array}{c} \Div \vc{V} \\ \vc{\omega} \times \vc{V} + \Grad s \end{array} \right]
\]
defined on the space $L^2(\Omega) \times L^2(\Omega;R^3)$. The operator $\mathcal{B}$ is skew symmetric, with the domain of definition
\[
\mathcal{D}[\mathcal{B}] = \left\{ [r, \vc{V}] \ \Big| \ r \in W^{1,2}(\Omega), \vc{V} \in L^2(\Omega;R^3), \Div \vc{V} \in L^2(\Omega)
\right\}.
\]

Next, we introduce the null space $\mathcal{N}[\mathcal{B}]$,
\bFormula{d3}
\mathcal{N}(\mathcal{B}) = \left\{ [q, \vc{v}] \ \Big| \ q = q(x_h), \ q \in W^{1,2}(R^2),
\ \vc{v} = [ \vc{v}_h (x_h), v_3],\ {\rm div}_h \vc{v}_h = 0 , \ \partial_{x_3} v_3 = 0, \right.
\eF
\[
\vc{\omega} \times \vc{v} + \Grad q = 0 \Big\}.
\]
We remark that $v_3 = 0$ as soon as $[q, \vc{v}] \in \mathcal{N}[\mathcal{B}]$ belongs to the symmetry class
(\ref{m2a}).

\subsubsection{Projection onto $\mathcal{N}(\mathcal{B})$}
\label{ker}

Given a couple of functions $[r, \vc{U}]$ defined in $\Omega$, we want to compute the orthogonal projection
\[
\mathcal{P} : L^2(\Omega) \times L^2(\Omega; R^3) \to \mathcal{N}(\mathcal{B}), \
\mathcal{P} [r, \vc{U}] = [q, \vc{v}].
\]
In addition, we assume the $[r, \vc{U}]$ belongs to the symmetry class (\ref{m2a}).

To begin, we project $[r, \vc{U}]$ onto the space of functions depending only on the horizontal variable $x_h$,
meaning we take
\[
\tilde r (x_h) = \frac{1}{2} \int_{-1}^1 r(x_h, x_3) \ {\rm d}x_3 = \int_0^1 r(x_h, x_3) \ {\rm d}x_3, \
\tilde \vc{U}_h (x_h) = \int_0^1 \vc{U}_h(x_h, x_3) \ {\rm d}x_3 ,\ \tilde U_3 = 0.
\]

Our goal is to minimize the functional
\[
[q, \vc{v}] \mapsto \| q - \tilde r \|^2_{L^2(R^2)} + \| \vc{v}_h - \tilde \vc{U}_h \|^2_{L^2(R^2; R^2)}
\]
under the constraint
\[
\vc{\omega} \times \vc{v} + \Grad q = 0.
\]

{
We have
\[
\left\| q - \tilde r \right\|^2_{L^2(R^2)} + \left\| \vc{v}_h -
\tilde \vc{U}_h \right\|^2_{L^2(R^2; R^2)}
\]
\[
=\left\| q - \tilde r \right\|^2_{L^2(R^2)}+
 \left\| (\vc{\omega}
\times \vc{v}) - (\vc{\omega} \times \tilde \vc{U})
\right\|^2_{L^2(R^2; R^2)}
\]
\[
= \left\| q - \tilde r \right\|^2_{L^2(R^2)}+ \left\| \nabla_h q +
(\vc{\omega} \times \tilde \vc{U}) \right\|^2_{L^2(R^2; R^2)}
\]
$$
=\int_{R^2}\Big(q^2+|\nabla_h q|^2- 2q\tilde r-2q \ {\rm
curl}_h\tilde \vc{U}_h\Big){\rm d} x_h+ \int_{R^2}\Big(\tilde r^2
+|\omega\times\tilde\vc U|^2\Big){\rm d}x_h.
$$
}
Thus the associated Euler-Lagrange equation for the minimization
problem give rise to \bFormula{d4} - \Delta_h q + q = \int_0^1
{\rm curl}_h \vc{U}_h \ {\rm d}x_3 + \int_0^1 r \ {\rm d}x_3, \eF
and \bFormula{d5} \vc{v} = [v_1,v_2],\ v_{1} = - \partial_{x_2}
q,\ v_{2} =  \partial_{x_1} q, \eF cf. the initial data
decomposition (\ref{m11}).

\subsection{Spectral analysis and dispersive estimates}

We employ the methods of \emph{Fourier analysis} in order to derive dispersive estimates for solutions of
the system (\ref{d1}), (\ref{d2}). Formally, the solutions of (\ref{d1}), (\ref{d2}) may be written in the form
\bFormula{d6}
[s,\vc{V}](t) = \exp \left( {\rm i} \frac{t}{\ep} [{\rm i}  \mathcal{B}] \right) [s_0, \vc{V}_0],
\eF
where ${\rm i}\mathcal{B}$ is a selfadjoint operator in $L^2(\Omega) \times L^2(\Omega;R^3)$.

Accordingly, we deduce that the solution operator generates a group of isometries in the $L^2$-norm,
specifically
\bFormula{d7}
\| s(t, \cdot) \|^2_{L^2(\Omega)} + \| \vc{V}(t, \cdot) \|^2_{L^2(\Omega;R^3)} =
\| s_0  \|^2_{L^2(\Omega)} + \| \vc{V}_0 \|^2_{L^2(\Omega;R^3)} \ \mbox{for any} \ t \in R.
\eF

Moreover, as the problem is linear, we obtain
\bFormula{d8}
\| s(t, \cdot) \|^2_{W^{m,2}(\Omega)} + \| \vc{V}(t, \cdot) \|^2_{W^{m,2}(\Omega;R^3)} =
\| s_0  \|^2_{W^{m,2}(\Omega)} + \| \vc{V}_0 \|^2_{W^{m,2}(\Omega;R^3)} \ \mbox{for any} \ t \in R,
\ m=0,1,\dots
\eF

\subsubsection{Fourier representation}

For each function $g \in L^2(\Omega)$, we introduce its Fourier representation
\[
\hat g (\xi, k), \ \xi = [\xi_1, \xi_2] \in R^2, \ k \in Z,
\]
where
\[
\hat g(\xi,k) = \frac{1}{\sqrt{2}} \int_{-1}^1 \int_{R^2} \exp
\left(-{\rm i} \xi \cdot x_h \right) g(x_h,x_3) \ {\rm d}x_h \exp
\left(- {\rm i} k  x_3 \right) \ {\rm d}x_3.
\]
We have
\[
g(x_h, x_3) = \sum_{k \in Z} \mathcal{F}^{-1}_{\xi \to x_h} \left[ \hat g(\xi,k ) \right] \exp \left({\rm i}k x_3 \right),
\]
where the symbol $\mathcal{F}_{x_h \to \xi}$ denotes the standard Fourier transform on $R^2$.

\subsubsection{Solutions in the Fourier variables}

The problem (\ref{d1}), (\ref{d2}) expressed in terms of the Fourier variables reads:
\bFormula{d9}
\ep \frac{{\rm d}}{{\rm d}t} \left[ \begin{array}{c} \hat s(t, \xi, k) \\ \hat {V}_1(t, \xi, k)
\\ \hat {V}_2(t, \xi, k) \\ \hat {V}_3(t, \xi, k)
 \end{array}
\right] + {\rm i} \left[ \begin{array}{cccc}
                   0 & \xi_1 & \xi_2 & k \\
                   \xi_1& 0 & {\rm i} & 0 \\
                   \xi_2& - {\rm i} & 0 & 0 \\
                    k& 0 & 0 & 0
                 \end{array}
 \right] \left[ \begin{array}{c} \hat s(t, \xi, k) \\ \hat {V}_1(t, \xi, k)
\\ \hat {V}_2(t, \xi, k) \\ \hat {V}_3(t, \xi, k)
 \end{array}
\right] = 0 , \left[ \begin{array}{c} \hat s(0, \xi, k) \\ \hat \vc{V}(0, \xi, k) \end{array}
\right] = \left[ \begin{array}{c} \hat s_0(\xi, k) \\ \hat \vc{V}_0(\xi, k) \end{array}
\right];
\eF
whence
\bFormula{d10}
[\hat s (t,\xi,k), \hat \vc{V}(t,\xi,k)] = \exp \left( - {\rm i} \frac{t}{\ep} \mathcal{A}(\xi,k) \right) [ \hat s_0 (\xi,k), \hat
\vc{V}_0 (\xi, k) ],
\eF
with the symmetric matrix
\[
\mathcal{A}(\xi,k) = \left[ \begin{array}{cccc}
                   0 & \xi_1 & \xi_2 & k \\
                   \xi_1& 0 & {\rm i} & 0 \\
                   \xi_2& - {\rm i} & 0 & 0 \\
                    k& 0 & 0 & 0
                 \end{array}
 \right].
\]

It is a routine matter to check that the symmetric matrix $\mathcal{A}(\xi,k)$ possesses four eigenvalues:
\bFormula{d11}
\begin{array}{c}
\lambda_{1} (\xi,k) = \left[ \frac{ 1 + | \xi |^2 + k^2 + \sqrt{ (1 + | \xi |^2 + k^2)^2 - 4 k^2 }}{2} \right]^{1/2} ,\ \lambda_2
(\xi, k)  = - \lambda_1 (\xi, k),
\\ \\
\lambda_{3} (\xi,k) = \left[ \frac{ 1 + | \xi |^2 + k^2 - \sqrt{ (1 + | \xi |^2 + k^2)^2 - 4 k^2 }}{2} \right]^{1/2},\ \lambda_4(\xi, k) = - \lambda_3 (\xi, k).
\end{array}
\eF
Note that $\lambda_3(\xi, 0) = \lambda_4(\xi, 0) = 0$, which corresponds to the non-trivial kernel of the operator $\mathcal{B}$ discussed in
Section \ref{ker}. Consequently, diagonalizing the matrix $\mathcal{A}$, we may rewrite (\ref{d10}) in the form
\bFormula{d12}
[\hat s(t,\xi,k), \hat \vc{V}(t,\xi, k)]
\eF
\[
= \mathcal{Q}^t (\xi, k) \exp \left( - \rm{i} \frac{t}{\ep}
\left[
         \begin{array}{cccc}
           \lambda_1(\xi,k) & 0 & 0 & 0 \\
            0& - \lambda_1 (\xi,k) & 0 & 0 \\
            0 & 0  & \lambda_3 (\xi, k)  & 0 \\
           0 & 0 & 0 & - \lambda_3 (\xi,k) \\
         \end{array}
       \right] \right) \mathcal{Q} (\xi, k) [ \hat s_0(\xi, k), \hat \vc{V}_0 (\xi, k) ]
\]
for a suitable matrix $\mathcal{Q}$.
\subsubsection{Decay estimates}
For each fixed $k$, the solution operators introduced in (\ref{d12}) may be viewed as
\[
\exp \left( - \rm{i} t \lambda_j (\xi, k) \right) = \exp \left( -
\rm{i} t \lambda_j (|\xi|, k) \right) {\approx} \exp \left( -
\rm{i} t \lambda_ j \left( \sqrt{ - \Delta_h } , k \right)
\right),\ j = 1,\dots,4.
\]
In particular, the eigenvalues are smooth functions of $|\xi| \approx \sqrt{ - \Delta_h }$ on the \emph{open} interval
$(0, \infty)$, and, as can be checked by direct computation, $\lambda_1(|\xi|, k)$ is strictly increasing in $|\xi|$
for any fixed $k \in Z$, while $\lambda_3(|\xi|, k)$ is strictly decreasing whenever $k \ne 0$. Consequently, we can use
the result of Guo, Peng, and Wang \cite[Theorem 1 (a)]{GuPeWa} to obtain the decay estimate
\bFormula{d13}
\left\| \exp \left( - \rm{i} t \lambda_j \left( \sqrt{- \Delta_h}, k \right) \right) {\Psi} \left(
\sqrt{-\Delta_h } \right) [v] \right\|_{L^\infty(R^2)} \leq |t|^{-1/2}
c(\Psi) \left\| {\Psi} \left(
\sqrt{-\Delta_h } \right) [v] \right\|_{L^1(R^2)},
\eF
\[
j = 1,2, \ k \in Z,\ j=3,4, \ k\ne 0,
\]
where,
\[
\Psi\left( \sqrt{-\Delta_h } \right)  [v] = \mathcal{F}^{-1}_{\xi
\to x_h} \left[ \Psi (|\xi|) \hat v(\xi) \right],\ \Psi \in \DC(0,\infty),
\]
is a frequency cut-off operator.

Knowing that $\exp \left( - \rm{i} t \lambda_j \left( \sqrt{-
\Delta_h}, k \right) \right)$ are also $L^2-$isometries, and using
the fact that  $\Psi (|\xi|)$ is $L^p$-multiplier, we conclude, by
interpolation, \bFormula{d14} \left\| \exp \left( - \rm{i} t
\lambda_j \left( \sqrt{- \Delta_h}, k \right) \right) {\Psi}
\left( \sqrt{-\Delta_h } \right) [v] \right\|_{L^p(R^2)} \leq
|t|^{\left(\frac{1}{p} - \frac{1}{2}\right)} c(\Psi,p) \| v
\|_{L^{p'}(R^2)}, \eF
\[
\mbox{for}\
p \geq 2, \ \frac{1}{p} + \frac{1}{p'} = 1, \mbox{and for}\
j = 1,2, \ k \in Z,\ j=3,4, \ k\ne 0,
\]

\subsection{Initial data decomposition}

In order to exploit the dispersive estimates derived in the preceding part, we have to find a suitable mollification of the initial data $\vr^{(1)}_0$ and
$\vu_0$. To this end we consider a family of smooth functions
\[
\psi_\delta \in \DC(0, \infty), 0 \leq \psi_\delta \leq 1, \ \psi_\delta \nearrow 1 \ \mbox{as}\ \delta \to 0,
\]
and
\[
\phi_\delta = \phi_\delta(x_h) \in \DC(R^2),\ 0 \leq \phi_\delta \leq 1, \ \phi_\delta \nearrow 1 \ \mbox{as}\ \delta \to 0.
\]

Finally, we regularize the data $\vr^{(1)}_0$, $\vu_0$ taking
\bFormula{d15} \left[ \vr^{(1)}_0 \right]_{\delta}(x_h,x_3) =
\frac 1{\sqrt 2}\sum_{|k|\le 1/\delta} \mathcal{F}^{-1}_{\xi \to
x_h} \left[ \psi_\delta (|\xi|) \widehat { \left( \vr_{0}^{(1)}
\phi_\delta \right) }(\xi,k) \right] \exp \left( -{\rm i} {k x_3}
\right),  \eF
 and, similarly, \bFormula{d16} \left[ u_{0,j}
\right]_{\delta} (x_h,x_3) =\frac 1{\sqrt 2}\sum_{|k|\le 1/\delta}
\mathcal{F}^{-1}_{\xi \to x_h} \left[ \psi_\delta (|\xi|) \widehat {
\left( u_{0,j} \phi_\delta   \right) }(\xi,k) \right] \exp \left(-
{\rm i} {k x_3} \right), \ j=1,2,3. \eF
 In other words,
we first multiply the data by a cut-off function to ensure
integrability and then perform a similar cut-off in the frequency
variable to ensure smoothness. We remark that
\begin{itemize}
\item the functions $\psi_\delta$ are obviously $L^p$ multipliers for any $1 < p < \infty$,
\item the orthogonal projection $\mathcal{P}$ onto the kernel $\mathcal{N}(\mathcal{B})$ commutes with the frequency cut-off represented by
$\psi_\delta$,
\item the operator ${\rm i}\mathcal{B}$ represented in the Fourier variables by the matrix $\mathcal{A}$ commutes with the frequency smoothing, in particular, the time evolution of the mollified data remains restricted to the domain of frequencies bounded above as well as below away from zero.
\end{itemize}

Finally, we write the initial data in the form
\bFormula{r2}
\left[ \vr^{(1)}_0 \right]_\delta = s_{0, \delta} + q_{0,\delta} ,\ \mbox{where} - \Delta_h q_{0, \delta} + q_{0, \delta} = \int_0^1 {\rm curl}_h \left[ [\vc{u}_0]_h \right]_\delta \ {\rm d}x_3 +
\int_0^1 \left[ \vr^{(1)}_0 \right]_{\delta} {\rm d}x_3
\eF
\bFormula{r3}
\left[ \vc{u}_0 \right]_{\delta} = \vc{V}_{0,\delta} + \vc{v}_{0, \delta} , \ \mbox{with}\ [v_{0,\delta}]_1 = - \partial_{x_2} q_{0, \delta},\ [v_{0,\delta}]_2 =  \partial_{x_1} q_{0, \delta},
\eF
and denote $s_{\ep, \delta}$, $\vc{V}_{\ep, \delta}$ the unique solution of the acoustic system
\bFormula{r6}
\ep \partial_t s_{\ep, \delta} + \Div \vc{V}_{\ep, \delta} = 0,
\eF
\bFormula{r7}
\ep \partial_t \vc{V}_{\ep, \delta} +  \vc{\omega} \times \vc{V}_{\ep, \delta} +  \Grad s_{\ep, \delta} = 0,
\eF
supplemented with the initial data
\bFormula{r8}
s_{\ep, \delta} (0, \cdot) = s_{0,\delta} ,\ \vc{V}_{\ep, \delta}(0, \cdot) = \vc{V}_{0,\delta}.
\eF

We have
\[
 [ \widehat{ s_{\ep,\delta}}(t,\xi,k),  \widehat{ \vc{V}_{\ep,\delta}}(t,\xi,
 k)]
 \]
 \[
= \mathcal{Q}^t (\xi, k) \exp \left( - \rm{i} \frac{t}{\ep} \left[
         \begin{array}{cccc}
           \lambda_1(\xi,k) & 0 & 0 & 0 \\
            0& - \lambda_1 (\xi,k) & 0 & 0 \\
            0 & 0  & \lambda_3 (\xi, k)  & 0 \\
           0 & 0 & 0 & - \lambda_3 (\xi,k) \\
         \end{array}
       \right] \right) \mathcal{Q} (\xi, k)
       \Big[ \widehat {s_{0,\delta}},\; \widehat {\vc{V}_{0,\delta}}
       \Big] (\xi, k),
\]
where, in agreement with the previous observations, {\bf (i)} each component of the vector
\[
\mathcal{Q} (\xi, k)
       \Big[ \widehat {s_{0,\delta}},\; \widehat {\vc{V}_{0,\delta}}
       \Big] (\xi, k)
\]
is of the form $\psi (|\xi|) a(\xi, k)$, with $\psi \in \DC(0, \infty)$ for
each fixed $k \in Z$, {\bf (ii)} as $[s_{0, \delta}, \vc{V}_{0, \delta}] \in \mathcal{N}[ \mathcal{B} ]^\perp$,
\[
\mathcal{Q} (\xi, 0)
       \Big[ \widehat {s_{0,\delta}},\; \widehat {\vc{V}_{0,\delta}}
       \Big] (\xi, 0) = \left[ \begin{array}{c} \psi_1(|\xi|) a_1(\xi) \\ \psi_2 (|\xi |) a_2(\xi) \\ 0 \\ 0 \end{array} \right].
\]

Consequently, interpolating the decay estimates (\ref{d14}) with (\ref{d8}) we may infer that
\bFormula{r9} s_{\ep, \delta} \to 0 \ \mbox{in}\ L^p(0,T; W^{m,
\infty}(\Omega)), \vc{V}_{\ep, \delta} \to 0 \ \mbox{in}\ L^p(0,T;
W^{m, \infty}(\Omega)) \ \mbox{as} \ \ep \to 0 \eF for any
\emph{fixed} $\delta > 0$, any $1 \leq p < \infty$, and
$m=0,1,\dots$.

\section{Convergence, part II}
\label{c}

We finish the proof of Theorem \ref{Tm1} by means of another application of the relative entropy inequality (\ref{r1}), this time for the choice
\bFormula{r4-}
r = r_{\ep, \delta} = 1 + \ep \Big( q_{\delta} + s_{\ep, \delta} \Big),\ \vc{U} = \vc{U}_{\ep, \delta} = \vc{v}_{\delta} + \vc{V}_{\ep, \delta},
\eF
where $[s_{\ep, \delta}, \vc{V}_{\ep, \delta}]$ are solutions of the acoustic system (\ref{r6} - \ref{r8}), the properties of which were discussed in the previous section, while $[q_{\delta}, \vc{v}_{\delta}]$ is the solution of the target problem
\bFormula{r4}
\vc{\omega} \times \vc{v}_{\delta} + \Grad q_{\delta} = 0
\eF
\bFormula{r5}
\partial_t \left( \Delta_h q_\delta -  q_\delta \right) + \nabla^\perp_h q_\delta \cdot \nabla_h
\left( \Delta_h q_\delta - q_\delta \right) = 0,\ q_\delta (0,
\cdot) = q_{0, \delta}, \eF cf. Proposition \ref{Pm1}.

\subsection{Initial data}

Going back to the relative entropy inequality (\ref{r1}) with the ``ansatz'' (\ref{r4-}) we get, in agreement with the hypotheses
(\ref{m11-} - \ref{m11a}),
\bFormula{r5a}
\mathcal{E}_\ep \left( \vr_{0,\ep}, \vu_{0,\ep} \ \Big| \ r(0,\cdot) , \vc{U}(0,\cdot) \right)
\eF
\[
=
\intO{ \frac{1}{2} \vr_{0,\ep} |\vu_{0,\ep} - [\vu_0]_\delta  |^2 }
\]
\[
+  \intO{ \left[ \frac{1}{\ep^2}
\left( H \left( 1 + \ep \vr^{(1)}_{0,\ep} \right) - \ep H'(1 + \ep [ \vr^{(1)}_0 ]_{\delta} ) \left((\vr^{(1)}_{0,\ep} - [\vr^{(1)}_0] _\delta
\right) - H (1 + \ep [ \vr^{(1)}_0 ]_{\delta} )\right) \right] }
\]
\[
\leq c \left( \| \vu_{0,\ep} - [ \vu_{0} ]_{\delta} \|^2_{L^2(\Omega;R^3)} +
\| \vr^{(1)}_{0,\ep} - [ \vr^{(1)}_{0} ]_{\delta} \|^2_{L^2(\Omega;R^3)} \right)
\]
\[
\to c \left( \| \vu_{0} - [ \vu_{0} ]_{\delta} \|^2_{L^2(\Omega;R^3)} +
\| \vr^{(1)}_{0} - [ \vr^{(1)}_{0} ]_{\delta} \|^2_{L^2(\Omega;R^3)} \right) \ \mbox{as} \ \ep \to 0.
\]

\subsection{Viscosity}

We write
\[
\left|
\int_0^\tau \intO{ \tn{S}_\ep (\Grad \vc{U}_{\ep, \delta} )  : \Grad (\vc{U}_{\ep,\delta} - \vu_\ep) } \ \dt \right|
\leq c_1 \int_0^\tau \intO{ \mu_\ep \left| \Grad (\vc{U}_{\ep,\delta} - \vu_\ep) \right| } \ \dt,
\]
where, by virtue of Korn's inequality,
\[
\int_0^\tau \intO{ \mu_\ep \left| \Grad (\vc{U}_{\ep,\delta} - \vu_\ep) \right| } \ \dt \leq
\frac{1}{2} \int_0^\tau \intO{ \left( \tn{S}_\ep(\Grad \vue) - \tn{S}_\ep (\Grad \vc{U}_{\ep, \delta} \right):
\Grad \left( \vc{U}_{\ep, \delta} - \vue \right) } \ \dt  + c_2 \mu_\ep.
\]

\subsection{Forcing term}

Furthermore, in accordance with the hypothesis (\ref{hypp}), the convergence established in (\ref{rr10}), and the
decay estimates (\ref{r9}), we get
\[
\int_0^\tau \intO{ \vre \Grad G \cdot ( \vc{U}_{\ep, \delta} - \vue ) } \ \dt \to \int_0^\tau \intO{  \Grad G \cdot ( \vc{v}_\delta - \vu ) } \ \dt = 0
\ \mbox{as} \ \ep \to 0,
\]
where we have used (\ref{rr14}), (\ref{r4}).

Summing up the previous estimates with
(\ref{r5a}), we can rewrite the relative entropy inequality in the form:

\bFormula{r10}
\intO{ \left[ \frac{1}{2} \vre |\vue - \vc{U}_{\ep, \delta}|^2 + \frac{1}{\ep^2} \left( H( \vre ) - H'(r_{\ep, \delta})(\vre - r_{\ep, \delta}) - H(r_{\ep, \delta}) \Big) \right) (\tau, \cdot) \right] }
\eF
\[
\leq
h_1 (\ep, \delta) + \int_0^\tau \intO{  \vre \left( \partial_t \vc{U}_{\ep, \delta} + \vue \cdot \Grad \vc{U}_{\ep, \delta} \right)
\cdot \left( \vc{U}_{\ep, \delta} - \vue \right) } \ \dt
\]
\[
 + \frac{1}{\ep} \int_0^\tau \intO{
\vre (\vc{\omega} \times \vue ) \cdot (\vc{U}_{\ep, \delta} - \vue) } \ \dt
\]
\[
+ \frac{1}{\ep^2} \int_0^\tau \intO{ \Big[ (r_{\ep, \delta} - \vre) \partial_t H'(r_{\ep, \delta}) + \Grad H'(r_{\ep, \delta}) \cdot (r_{\ep, \delta} \vc{U}_{\ep, \delta} - \vre \vue ) \Big] } \ \dt
\]
\[
- \frac{1}{\ep^2} \int_0^\tau \intO{ \Div \vc{U}_{\ep, \delta} \Big( p(\vre) - p(r_{\ep, \delta}) \Big) } \ \dt .
\]
Here and hereafter, we use the symbol $h_i, i=1,2,\dots$ to denote a function of $\ep$, $\delta$ enjoying the following properties
\[
h_i(\ep, \delta) \to \tilde h_i(\delta) \ \mbox{for}\ \ep \to 0, \ \mbox{with}\ \tilde h_i(\delta) \to 0
\ \mbox{for}\ \delta \to 0.
\]

\subsection{Estimating the remaining terms}

To begin, let us recall our convention that 
\[
\Ov{\vr} = p'(\Ov{\vr}) = H''(\Ov{\vr}) = 1.
\]
Furthermore, in the following discussion, 
we make a systematic use of the dispersive decay estimates (\ref{r9}).

\bigskip

{\bf Step 1:}

\medskip

We have
\bFormula{r11}
\Big[ (r_{\ep,\delta} - \vre) \partial_t H'(r_{\ep, \delta}) + \Grad H'(r_{\ep, \delta}) \cdot (r_{\ep, \delta} \vc{U}_{\ep, \delta} - \vre \vue ) \Big] - \Div \vc{U}_{\ep, \delta} \Big( p(\vre) - p(r_{\ep, \delta}) \Big)
\eF
\[
= \Big[ p(r_{\ep, \delta})  - p'(r_{\ep, \delta}) (r_{\ep, \delta} - \vre) - p(\vre) \Big] \Div \vc{U}_{\ep, \delta}
\]
\[
+ (r_{\ep, \delta} - \vre) H''(r_{\ep, \delta}) \Big[ \partial_t r_{\ep, \delta} + \Div (r_{\ep, \delta} \vc{U}_{\ep, \delta}) \Big] + \vre \Grad H'(r_{\ep, \delta}) \cdot ( \vc{U}_{\ep, \delta} - \vue),
\]
where, in accordance with (\ref{r6}), (\ref{r4}),
\bFormula{r12}
\partial_t r_{\ep, \delta} + \Div (r_{\ep, \delta} \vc{U}_{\ep, \delta}) = \ep \partial_t q_\delta + \ep \partial_t s_{\ep, \delta} + \Div(r_{\ep, \delta} (\vc{v}_\delta + \vc{V}_{\ep, \delta} )) =
\ep \partial_t q_\delta + \ep \Div \Big[ (q_\delta + s_{\ep, \delta} ) \vc{U}_{\ep, \delta} \Big].
\eF

Furthermore, we get
\bFormula{r13}
\vre \Grad H'(r) \cdot ( \vc{U}_{\ep, \delta} - \vue) = \vre \Grad \Big[ H'(r_{\ep, \delta}) - H''(\Ov{\vr})(r_{\ep, \delta} - \Ov{\vr}) - H'(\Ov{\vr}) \Big]\cdot ( \vc{U}_{\ep, \delta} - \vue)
\eF
\[
+ \ep \vre \frac{p'(\Ov{\vr})}{\Ov{\vr}} \Grad q_\delta \cdot (\vc{U}_{\ep, \delta} - \vue) +
\ep \vre \frac{p'(\Ov{\vr})}{\Ov{\vr}} \Grad s_{\ep, \delta} \cdot (\vc{U}_{\ep, \delta} - \vue)
\]
\[
= \vre \Grad \Big[ H'(r_{\ep, \delta}) - H''(1)(r_{\ep, \delta} - 1) - H'(1) \Big]\cdot ( \vc{U}_{\ep, \delta} - \vue)
\]
\[
+ \ep \vre \frac{p'(\Ov{\vr})}{\Ov{\vr}} \Grad q_\delta \cdot (\vc{U}_{\ep, \delta} - \vue) - \ep^2 \vre (\vc{U}_{\ep, \delta} - \vue) \cdot \partial_t \vc{V}_{\ep, \delta}
\]
\[
- \ep \vre (\vc{U}_{\ep, \delta} - \vue) \cdot (\vc{\omega} \times \vc{V}_{\ep, \delta}  )
\]

Consequently, after a straightforward manipulation, the inequality (\ref{r10}) can be rewritten as follows:
\[
\intO{ \left[ \frac{1}{2} \vre |\vue - \vc{U}_{\ep, \delta}|^2 + \frac{1}{\ep^2} \Big( H(\vre) - H'(r_{\ep, \delta})(\vre - r_{\ep, \delta}) - H(r_{\ep, \delta}) \Big)(\tau, \cdot) \right] }
\]
\[
\leq
h_1(\ep, \delta) + \int_0^\tau \intO{  \vre \left( \partial_t \vc{v}_\delta + \vue \cdot \Grad \vc{U}_{\ep, \delta} \right) \cdot \left(
\vc{U}_{\ep, \delta} - \vue \right) } \ \dt
\]
\[
 + \frac{1}{\ep} \int_0^\tau \intO{
\vre (\vc{\omega} \times \vc{v}_\delta ) \cdot (\vc{U}_{\ep, \delta} - \vue) } \ \dt
\]
\[
+ \frac{1}{\ep^2} \int_0^\tau \intO{ \Big[ p(r_{\ep, \delta})  - p'(r_{\ep, \delta}) (r_{\ep, \delta} - \vre) - p(\vre) \Big] \Div \vc{U}_{\ep, \delta}  } \ \dt
\]
\[
+ \frac{1}{\ep} \int_0^\tau \intO{ (r_{\ep, \delta} - \vre) H''(r_{\ep, \delta}) \left[ \partial_t q_\delta + \Div \Big( (q_\delta + s_{\ep,
\delta} ) \vc{U}_{\ep, \delta} \Big) \right] } \ \dt
\]
\[
+ \frac{1}{\ep^2} \int_0^\tau \intO{\vre \Grad \Big[ H'(r_{\ep, \delta}) - H''(1)(r_{\ep, \delta} - 1) - H'(1) \Big]\cdot ( \vc{U}_{\ep, \delta} - \vue)} \ \dt
\]
\[
+ \frac{1}{\ep} \int_0^\tau \intO{ \vre \frac{p'(\Ov{\vr})}{\Ov{\vr}} \Grad q_\delta \cdot (\vc{U}_{\ep, \delta} - \vue) } \ \dt
\]

Finally, we use the relation (\ref{r4}) to conclude that
\bFormula{r14}
\intO{ \left[ \frac{1}{2} \vre |\vue - \vc{U}_{\ep, \delta}|^2 + \frac{1}{\ep^2} \Big( H(\vre) - H'(r_{\ep, \delta})(\vre - r_{\ep, \delta}) - H(r_{\ep, \delta}) \Big)(\tau, \cdot) \right] }
\eF
\[
\leq
h_1(\ep, \delta) + \int_0^\tau \intO{  \vre \left( \partial_t \vc{v}_\delta + \vue \cdot \Grad \vc{U}_{\ep, \delta} \right) \cdot \left(
\vc{U}_{\ep, \delta} - \vue \right) } \ \dt
\]
\[
+ \frac{1}{\ep^2} \int_0^\tau \intO{ \Big[ p(r_{\ep, \delta})  - p'(r_{\ep, \delta}) (r_{\ep, \delta} - \vre) - p(\vre) \Big] \Div \vc{U}_{\ep, \delta}  } \ \dt
\]
\[
+ \frac{1}{\ep} \int_0^\tau \intO{ (r_{\ep, \delta} - \vre) H''(r_{\ep, \delta}) \left[ \partial_t q_\delta + \Div \Big( (q_\delta + s_{\ep,
\delta} ) \vc{U}_{\ep, \delta} \Big) \right] } \ \dt
\]
\[
+ \frac{1}{\ep^2} \int_0^\tau \intO{\vre \Grad \Big[ H'(r_{\ep, \delta}) - H''(1)(r_{\ep, \delta} - 1) - H'(1) \Big]\cdot ( \vc{U}_{\ep, \delta} - \vue)} \ \dt.
\]

\medskip

\noindent {\bf Step 2:}

\medskip

We rewrite the integral
\[
\int_0^\tau \intO{  \vre \left( \partial_t \vc{v}_\delta + \vue \cdot \Grad \vc{U}_{\ep, \delta} \right) \cdot \left( \vc{U}_{\ep, \delta} -
\vue \right) } \ \dt =
\int_0^\tau \intO{  \vre \left( \partial_t \vc{v}_\delta + \vc{v}_\delta \cdot \Grad \vc{v}_\delta \right) \cdot \left( \vc{U}_{\ep, \delta} - \vue \right) } \ \dt
\]
\[
+ \int_0^\tau \intO{ \vre \vc{v}_\delta \cdot \Grad \vc{V}_{\ep, \delta} \cdot (\vc{U}_{\ep, \delta} - \vue) } \ \dt +
\int_0^\tau \intO{ \vre \vc{V}_{\ep, \delta} \cdot \Grad \vc{U}_{\ep, \delta} \cdot (\vc{U}_{\ep, \delta} - \vue) } \ \dt
\]
\[
- \int_0^\tau \intO{  \vre (\vc{U}_{\ep, \delta} - \vue) \cdot \Grad \vc{U}_{\ep,\delta} \cdot (\vc{U}_{\ep, \delta} - \vue) } \ \dt,
\]
where, by virtue of the uniform bounds (\ref{rr1} - \ref{rr4}), combined with the dispersive estimates (\ref{r9}),
\[
\left| \int_0^\tau \intO{ \vre \vc{v}_\delta \cdot \Grad \vc{V}_{\ep, \delta} \cdot (\vc{U}_{\ep, \delta} - \vue) } \ \dt +
\int_0^\tau \intO{ \vre \vc{V}_{\ep, \delta} \cdot \Grad \vc{U}_{\ep, \delta} \cdot (\vc{U}_{\ep, \delta} - \vue) } \ \dt \right|
\to 0 \ \mbox{as} \ \ep \to 0
\]
uniformly in $\tau \in [0,T]$.

Consequently, we may infer that
\bFormula{r15}
\mathcal{E}_\ep \left( \vre , \vue \Big| r_{\ep, \delta}, \vc{U}_{\ep, \delta} \right) \equiv
\intO{ \left[ \frac{1}{2} \vre |\vue - \vc{U}_{\ep, \delta} |^2 + \frac{1}{\ep^2} \Big( H(\vre) - H'(r_{\ep, \delta})(\vre - r_{\ep, \delta}) - H(r_{\ep, \delta}) \Big)(\tau, \cdot) \right] }
\eF
\[
\leq
h_2 (\ep, \delta) + {\rm sup}_{t \in [0,T]} \| \vc{v}_\delta (t, \cdot) \|_{W^{1, \infty}(R^2; R^2)} \int_0^\tau  \mathcal{E}_\ep \left( \vre , \vue \Big| r_{\ep, \delta}, \vc{U}_{\ep, \delta} \right) (t, \cdot) \ {\rm d}t
\]
\[
+ \int_0^\tau \intO{  \vre \left( \partial_t \vc{v}_\delta + \vc{v}_\delta \cdot \Grad \vc{v}_\delta \right) \cdot \left( \vc{U}_{\ep, \delta} - \vue \right) } \ \dt
\]
\[
+ \frac{1}{\ep^2} \int_0^\tau \intO{ \Big[ p(r_{\ep, \delta})  - p'(r_{\ep, \delta}) (r_{\ep, \delta} - \vre) - p(\vre) \Big] \Div \vc{U}_{\ep, \delta}  } \ \dt
\]
\[
+ \frac{1}{\ep} \int_0^\tau \intO{ (r_{\ep, \delta} - \vre) H''(r_{\ep, \delta}) \left[ \partial_t q_\delta + \Div \Big( (q_{\delta} +
s_{\ep, \delta}) \vc{U}_{\ep, \delta} \Big) \right] } \ \dt
\]
\[
+ \frac{1}{\ep^2} \int_0^\tau \intO{\vre \Grad \Big[ H'(r_{\ep, \delta}) - H''(1)(r_{\ep, \delta} - 1) - H'(1) \Big]\cdot ( \vc{U}_{\ep, \delta} - \vue)} \ \dt.
\]

\medskip

\noindent
{\bf Step 3:}

\medskip

In view of the uniform bounds (\ref{rr2} - \ref{rr4}), we have
\[
{\rm ess} \sup_{t \in (0,T)} \frac{1}{\ep^2} \left\| \left[ p(r_{\ep, \delta})  - p'(r_{\ep, \delta}) (r_{\ep, \delta} - \vre) - p(\vre)
 \right] (t, \cdot) \right\|_{L^1(\Omega)}
\leq c;
\]
whence, by virtue of the dispersive estimates (\ref{r9}),
\bFormula{r16}
\frac{1}{\ep^2} \int_0^\tau \intO{ \Big[ p(r_{\ep, \delta})  - p'(r_{\ep, \delta}) (r_{\ep, \delta} - \vre) - p(\vre) \Big] \Div \vc{U}_{\ep, \delta}  } \ \dt
\eF
\[
= \frac{1}{\ep^2} \int_0^\tau \intO{ \Big[ p(r_{\ep, \delta})  - p'(r_{\ep, \delta}) (r_{\ep, \delta} - \vre) - p(\vre) \Big] \Div \vc{V}_{\ep, \delta}  } \ \dt
\to 0 \ \mbox{as}\ \ep \to 0,
\]
uniformly for $\tau \in [0,T]$.

Similarly, we obtain
\[
\frac{1}{\ep^2} \Grad \Big[ H'(r_{\ep, \delta}) - H''(1)(r_{\ep, \delta} - 1) - H'(1) \Big]
\]
\[
= \frac{1}{\ep} \left( H''(r_{\ep, \delta}) - H''(1) \right)
\Grad ( q_\delta + s_{\ep, \delta} ) \to H'''(1) q_\delta \Grad q_\delta \ \mbox{in}\ L^\infty(0,T; L^2 \cap L^\infty(\Omega; R^3))
\ \mbox{as}\ \ep \to 0.
\]
Thus, using (\ref{rr10}), (\ref{rr11}) we conclude that
\bFormula{r17}
\frac{1}{\ep^2} \int_0^\tau \intO{\vre \Grad \Big[ H'(r_{\ep, \delta}) - H''(1)(r_{\ep, \delta} - 1) - H'(1) \Big]\cdot ( \vc{U}_{\ep, \delta} - \vue)} \ \dt
\to 0 \ \mbox{as}\ \ep \to 0
\eF
uniformly in $\tau \in [0,T]$.

In view of (\ref{r16}), (\ref{r17}), the relative entropy inequality (\ref{r15}) reduces to
\bFormula{r18}
\mathcal{E}_\ep \left( \vre , \vue \Big| r_{\ep, \delta}, \vc{U}_{\ep, \delta} \right) \equiv
\intO{ \left[ \frac{1}{2} \vre |\vue - \vc{U}_{\ep, \delta} |^2 + \frac{1}{\ep^2} \Big( H(\vre) - H'(r_{\ep, \delta})(\vre - r_{\ep, \delta}) - H(r_{\ep, \delta}) \Big)(\tau, \cdot) \right] }
\eF
\[
\leq
h_3 (\ep, \delta) + {\rm sup}_{t \in [0,T]} \| \vc{v}_\delta (t, \cdot) \|_{W^{1, \infty}(R^2; R^2)} \int_0^\tau  \mathcal{E}_\ep \left( \vre , \vue \Big| r_{\ep, \delta}, \vc{U}_{\ep, \delta} \right) (t, \cdot) \ {\rm d}t
\]
\[
+ \int_0^\tau \intO{  \vre \left( \partial_t \vc{v}_\delta + \vc{v}_\delta \cdot \Grad \vc{v}_\delta \right) \cdot \left( \vc{U}_{\ep, \delta} - \vue \right) } \ \dt
\]
\[
+ \frac{1}{\ep} \int_0^\tau \intO{ (r_{\ep, \delta} - \vre) H''(r_{\ep, \delta}) \left[ \partial_t q_\delta + \Div \Big( (q_{\delta} +
s_{\ep, \delta}) \vc{U}_{\ep, \delta} \Big) \right] } \ \dt.
\]

\medskip

\noindent
{\bf Step 4:}

\medskip

Similarly to the above, we deduce from (\ref{r18}) that
\bFormula{r19}
\mathcal{E}_\ep \left( \vre , \vue \Big| r_{\ep, \delta}, \vc{U}_{\ep, \delta} \right) \equiv
\intO{ \left[ \frac{1}{2} \vre |\vue - \vc{U}_{\ep, \delta} |^2 + \frac{1}{\ep^2} \Big( H(\vre) - H'(r_{\ep, \delta})(\vre - r_{\ep, \delta}) - H(r_{\ep, \delta}) \Big)(\tau, \cdot) \right] }
\eF
\[
\leq
h_4 (\ep, \delta) + {\rm sup}_{t \in [0,T]} \| \vc{v}_\delta (t, \cdot) \|_{W^{1, \infty}(R^2; R^2)} \int_0^\tau  \mathcal{E}_\ep \left( \vre , \vue \Big| r_{\ep, \delta}, \vc{U}_{\ep, \delta} \right) (t, \cdot) \ {\rm d}t
\]
\[
+ \int_0^\tau \intO{  \left( \partial_t \vc{v}_\delta + \vc{v}_\delta \cdot \Grad \vc{v}_\delta \right) \cdot \left( \vc{v}_\delta - \vu \right) } \ \dt
+ \int_0^\tau \intO{ (q_\delta - \vr^{(1)}) \left[ \partial_t q_\delta + \Div \Big( q_{\delta}
\vc{v}_\delta \Big) \right] } \ \dt,
\]
where
\[
\intO{  \left( \partial_t \vc{v}_\delta + \vc{v}_\delta \cdot \Grad \vc{v}_\delta \right) \cdot \left( \vc{v}_\delta - \vu \right) }
+ \intO{ (q_\delta - \vr^{(1)}) \left[ \partial_t q_\delta + \Div \Big( q_{\delta}
\vc{v}_\delta \Big) \right] }
\]
\[
= \frac{1}{2} \frac{{\rm d}}{{\rm d}t} \intO{ \left( |\vc{v}_\delta|^2 + {q}^2_{\delta} \right) } -
\intO{ \left( \partial_t \vc{v}_\delta \cdot \vu + \partial_t q_\delta \vr^{(1)} \right) }
\]
\[
- \intO{ \left( \vc{v}_\delta \cdot \Grad \vc{v}_\delta \cdot \vu + \vr^{(1)} \Div (q_\delta \vc{v}_\delta) \right) }.
\]

Furthermore, we have
\[
\Div (q_\delta \vc{v}_\delta) = \Grad q_\delta \cdot \vc{v}_\delta = -\Grad q_\delta \cdot \Grad^t q_\delta = 0,
\]
and, by virtue of (\ref{m8}), (\ref{rr12}), 
\[
-
\intO{ \left( \partial_t \vc{v}_\delta \cdot \vu + \partial_t q_\delta \vr^{(1)} \right) } =
\intO{ \partial_t \left( \Delta_h q_\delta - q_\delta \right) \vr^{(1)} } = - \intO{ \vr^{(1)} \nabla_h ( \Delta_h q_\delta ) \cdot \vc{v}_\delta }
\]
\[
=
- \intO{ (\vc{\omega} \times \vu) \cdot \vc{v}_\delta \Delta_h q_\delta }.
\]

On the other hand, a routine manipulation gives rise to
\[
\vc{v}_\delta \cdot \Grad \vc{v}_\delta \cdot \vu +
(\vc{\omega} \times \vu) \cdot \vc{v}_\delta \Delta_h q_\delta =  \vu \cdot \nabla_h |\vc{v}_\delta|^2;
\]
whence
\[
\intO{  \left( \partial_t \vc{v}_\delta + \vc{v}_\delta \cdot \Grad \vc{v}_\delta \right) \cdot \left( \vc{v}_\delta - \vu \right) }
+ \intO{ (q_\delta - \vr^{(1)}) \left[ \partial_t q_\delta + \Div \Big( q_{\delta}
\vc{v}_\delta \Big) \right] }
\]
\[
= \frac{1}{2} \frac{{\rm d}}{{\rm d}t} \intO{ \left( |\vc{v}_\delta|^2 + {q}^2_{\delta} \right) }=
\frac{1}{2} \frac{{\rm d}}{{\rm d}t} \intO{ \left( |\nabla_h q_\delta |^2 + {q}^2_{\delta} \right) }
\]

Finally, multiplying the equation (\ref{r5}) on $q_\delta$ we obtain
\[
\frac{1}{2} \frac{{\rm d}}{{\rm d}t} \intO{ \left( |\nabla_h q_\delta |^2 + q^2_{\delta} \right) } =
\intO{ \vc{v} \cdot \nabla_h (\Delta_h q_\delta ) q_\delta }
= - \intO{ \vc{v} \cdot \nabla_h  q_\delta \Delta_h q_\delta } = 0.
\]

In view of the previous discussion, the relation (\ref{r19}) reduces to
\bFormula{r20}
\mathcal{E}_\ep \left( \vre , \vue \Big| r_{\ep, \delta}, \vc{U}_{\ep, \delta} \right) \equiv
\intO{ \left[ \frac{1}{2} \vre |\vue - \vc{U}_{\ep, \delta} |^2 + \frac{1}{\ep^2} \Big( H(\vre) - H'(r_{\ep, \delta})(\vre - r_{\ep, \delta}) - H(r_{\ep, \delta}) \Big)(\tau, \cdot) \right] }
\eF
\[
\leq
h_4 (\ep, \delta) + {\rm sup}_{t \in [0,T]} \| \vc{v}_\delta (t, \cdot) \|_{W^{1, \infty}(R^2; R^2)} \int_0^\tau  \mathcal{E}_\ep \left( \vre , \vue \Big| r_{\ep, \delta}, \vc{U}_{\ep, \delta} \right) (t, \cdot) \ {\rm d}t,
\]
with
\[
h_4(\ep, \delta) \to \tilde h_4(\delta) \ \mbox{for}\ \ep \to 0, \ \mbox{with}\ \tilde h_4(\delta) \to 0
\ \mbox{for}\ \delta \to 0.
\]

Passing to the limit, first for $\ep \to 0$ and then for $\delta \to 0$, in (\ref{r20}), we complete the proof of Theorem \ref{Tm1}.

\def\cprime{$'$} \def\ocirc#1{\ifmmode\setbox0=\hbox{$#1$}\dimen0=\ht0
  \advance\dimen0 by1pt\rlap{\hbox to\wd0{\hss\raise\dimen0
  \hbox{\hskip.2em$\scriptscriptstyle\circ$}\hss}}#1\else {\accent"17 #1}\fi}


\end{document}